\documentclass[twoside,12pt]{amsart}
\usepackage{amsfonts}

\usepackage{graphicx}
\usepackage{amscd}
\usepackage{amsmath}
\begin{document}


\newtheorem{Thm}{Theorem}
\newtheorem{Ax}{Axiom}
\newtheorem{Prop}{Proposition}
\newtheorem{Cor}[Prop]{Corollary}
\newtheorem{Main}{}
\renewcommand{\theMain}{}
\newtheorem{Lem}[Prop]{Lemma}
\newtheorem{Fact}{Fact}
\renewcommand{\theFact}{}
\newtheorem{Claim}{Claim}
\renewcommand{\theClaim}{}

\newtheorem{Def}{Definition}
\newtheorem{rmk}{Remark}
\newenvironment{Rmk}{\begin{rmk}\em}{\end{rmk}}
\newtheorem{exm}{Example}
\newenvironment{Exm}{\begin{exm}\em}{\end{exm}}

\newtheorem{prf}{Proof}
\renewcommand{\theprf}{}
\newenvironment{Prf}{\begin{prf}\em}{\qed\end{prf}}
\newtheorem{prff}{}
\renewcommand{\theprff}{}
\newenvironment{Prff}{\begin{prff}\em}{\qed\end{prff}}



\newcommand{\YES}[1]{#1}
\newcommand{\NOT}[1]{}

\newcommand{\cA}{\mathcal{A}}
\newcommand{\cB}{\mathcal{B}}
\newcommand{\cC}{\mathcal{C}}
\newcommand{\cD}{\mathcal{D}}
\newcommand{\cE}{\mathcal{E}}
\newcommand{\cF}{\mathcal{F}}
\newcommand{\cG}{\mathcal{G}}
\newcommand{\cH}{\mathcal{H}}
\newcommand{\cI}{\mathcal{I}}
\newcommand{\cJ}{\mathcal{J}}
\newcommand{\cK}{\mathcal{K}}
\newcommand{\cL}{\mathcal{L}}
\newcommand{\cM}{\mathcal{M}}
\newcommand{\cN}{\mathcal{N}}
\newcommand{\cO}{\mathcal{O}}
\newcommand{\cP}{\mathcal{P}}
\newcommand{\cQ}{\mathcal{Q}}
\newcommand{\cR}{\mathcal{R}}
\newcommand{\cS}{\mathcal{S}}
\newcommand{\cT}{\mathcal{T}}
\newcommand{\cU}{\mathcal{U}}
\newcommand{\cV}{\mathcal{V}}
\newcommand{\cW}{\mathcal{W}}
\newcommand{\cX}{\mathcal{X}}
\newcommand{\cY}{\mathcal{Y}}
\newcommand{\cZ}{\mathcal{Z}}

\newcommand{\bbb}[1]{{\mathbb{#1}}}

\newcommand{\bN}{\bbb{N}}
\newcommand{\bZ}{\bbb{Z}}
\newcommand{\bR}{\bbb{R}}
\newcommand{\bC}{\bbb{C}}
\newcommand{\bQ}{\bbb{Q}}
\newcommand{\bT}{\bbb{T}}

\newcommand{\noind}[1]{{\setlength{\parindent}{0cm} #1}}
\newcommand{\parsk}{\par\medskip}
\newcommand{\pa}{\par\medskip}

\newcommand{\varend}{\end{document}}

\newcommand{\LP}{\left(}  \newcommand{\RP}{\right)}
\newcommand{\LQ}{\left[}  \newcommand{\RQ}{\right]}
\newcommand{\LB}{\left\{} \newcommand{\RB}{\right\}}
\newcommand{\LA}{\left\langle} \newcommand{\RA}{\right\rangle}
\newcommand{\LS}{\left|}  \newcommand{\RS}{\right|}
\newcommand{\LN}{\left\|} \newcommand{\RN}{\right\|}
\newcommand{\bLP}{\bigl(}  \newcommand{\bRP}{\bigr)}

\newcommand{\ts}{{\otimes}} \newcommand{\TS}{{\bigotimes}}

\newcommand{\beware}{~$\Psi\!\Psi\!\Psi$~}

\newcommand{\EM}{\em\bf}

\newcommand{\BE}{\begin{equation}}  \newcommand{\EE}{\end{equation}}
\newcommand{\BR}{\begin{eqnarray*}}  \newcommand{\ER}{\end{eqnarray*}}
\newcommand{\BER}{$$\begin{array}}  \newcommand{\EER}{\end{array}$$}
\newcommand{\head}[1]{
             \par\bigskip\noind{{\large\bf#1}\nopagebreak\par}\medskip}
\newcommand{\chead}[1]{\begin{center}
              \par\bigskip{\large\bf#1}\nopagebreak\par\medskip
              \end{center}}

\newcommand{\fillitem}{\L{\embox{}}\vspace{-.6cm}}

\newcommand{\toprove}[1]{{\buildrel\mbox{?}\over#1}}
\newcommand{\DEF}{{\buildrel\mbox{\scriptsize\,\,def\,\,}\over=}}
\newcommand{\RRe}{{\mbox{Re}\,}}
\newcommand{\IIm}{{\mbox{Im}\,}}
\newcommand{\BAR}[1]{\overline{#1}}
\newcommand{\HAT}[1]{{\buildrel{\,\,\wedge}\over{#1}}\!}
\newcommand{\TIL}[1]{{\buildrel{\,\sim\,}\over{#1}}}

\newcommand{\I}{\infty}
\newcommand{\es}{\emptyset}
\newcommand{\sm}{\setminus} \newcommand{\st}{\subset}
\newcommand{\eps}{\varepsilon}
\newcommand{\al}{\alpha} \newcommand{\be}{\beta}
\newcommand{\ga}{\gamma} \newcommand{\de}{\delta}
\newcommand{\DE}{\Delta}
\newcommand{\OM}{\Omega} \newcommand{\om}{\omega}
\newcommand{\la}{\lambda}
\newcommand{\half}{\frac{1}{2}}
\newcommand{\NI}{{n\to\infty}}
\newcommand{\ang}{{<\!\!\!)}}
\newcommand{\arc}[1]{\breve{#1}}

\newenvironment{txeqn}[1]{\begin{equation}\label{#1}\begin{array}{l}}{
                         \end{array}\end{equation}}


\newcommand{\card}{{\text{card}\,}}
\newcommand{\lr}{{lim-rim}}
\newcommand{\Lr}{{Lim-Rim}}
\newcommand{\T}{^*\!\!}
\newcommand{\U}{^*\!}
\newcommand{\NS}{non-standard}
\newcommand{\NNS}{Non-Standard}
\newcommand{\NA}{Non-Standard Analysis}

\title{An Approach to Non-Standard Analysis}
\author{Eliahu Levy}
\address{Department of Mathematics,
Technion -- Israel Institute of Technology,
Haifa 32000, Israel}
\email{eliahu@techunix.technion.ac.il}


\date{}



\begin{abstract}
This note has two principal aims: to portray an essence of \NA\ as a
particular structure (which we call \lr), noting its interplay with the
notion of ultrapower, and to present a construction of \NA, viewed as a matter
of mathematics, where the set $\U A$ of non-standard elements of a set $A$ --
``the adjunction of all possible limits'' is a ``good'' kind of \lr\ which
plays a role analogous to that of the algebraic closure of a field -- ``the
adjunction of all roots of polynomials''.
In the same spirit as with algebraic closures, one has uniqueness
up to isomorphism, and also universality and homogeneity, provided one has
enough General Continuum Hypothesis.
The cardinality of $\U A$ will be something like $2^{2^{|A|}}$ -- the same as
that of the set of ultrafilters in $A$, and one has a high degree of
saturation.
\end{abstract}

\maketitle

\setcounter{section}{-1}

\section{Notations}

Denote the cardinality of a set $A$ by $|A|$.
For a cardinality $m$, denote by $m^+$ the successor cardinality.

\section{Features of the Conventional Approach to \NA}
Non-standard analysis was constructed by its founder, Abraham Robinson
(\cite{Robinson}) as a non-standard model of part of Set Theory. In his and
following treatments (see several treatises in the bibliography)%
\footnote{In another direction, beginning with E.~Nelson \cite{Nelson}, one
constructs a special Set Theory for \NA.}
the idea is to start with a ``big'' structure -- a set $V$,
with the restriction to it of the membership relation $\in$, which will be
big enough so that usual mathematics could be done in it (should contain,
e.g.\ $\bN$ (the natural numbers), $\bR$ (the real numbers), sets of them,
sets of these etc.) $V$ is to be provided with a ``mirror'' ``non-standard''
structure $^*V$ with a relation $\T\in$, viewed as an ``interpretation'' of
the $\in$ of $V$ (Also, ``interpret'' ``$=$'' in $V$ by ``$=$'' in $^*V$),
so that

\begin{enumerate}
\item
$V$ is embedded in $^*V$ -- identified with a subset of $^*V$. The members of
$V$ (thus considered as members of $^*V$) are referred to as {\EM standard}
(yet by ``\NS'' elements we mean any members of $^*V$).

\item
A ``Transfer Principle'' holds: any (in some constructions one has to say:
``bounded'')
first-order logical expression with members of $V$ as constants
which is a true sentence when quantifiers are interpreted relative to $V$
will remain a true sentence if (the constants -- members of $V$ -- are replaced by
their identified images in $^*V$ and) $\in$ is replaced by $\T\in$ and
the quantifiers are interpreted relative to $^*V$.
In this sense $^*V$ is a model to the true in $V$ sentences of the
first-order theory with $\in$ and all elements of $V$ as constants.

\item
By Transfer, for any standard set $A$, a standard $a\in V$ will be
an $\T\in$-member of it iff it is a usual member, but in $^*V$\,\,$A$ may
have other, {\em \NS} members. The set of all elements of $^*V$ which
$\T\in$-belong to $A$ is denoted by $\U A$. Again by Transfer, $^*(A\times B)$
is canonically identified with $\U A\times\,^*B$, thus for relations, functions
etc.\ we have the $^*$-relation, $^*$-function etc.\ among members of $^*V$
(which for standard elements coincides with the original relation, function
etc.) When there is no danger of misunderstanding, one often uses the same
symbol (e.g.\ $+$, $<$) for the $^*$-relation etc.\ as for the original.

\item
$^*V$ is indeed ``\NS'' -- it contains elements not in $V$,
and one has features, usually some kind of ``saturation'' (see below),
which will ensure, say, that there is an $^*n\in{}^*V$ which $\T\in$-belongs
to the set $\bN$ of natural numbers, but such that for every standard
$m\in N$,\,\,$m\;\T<{}^*n$ (i.e.\ $(m,{}^*n)$ $\T\in$-belongs to the graph of
$<$). Such $^*n$'s are naturally referred to as
``infinite'' or ``unlimited''%
\footnote{In fact, for $\bN$ the existence of ``infinite'' members needs only
the existence of some \NS\ members -- one easily proves that all
$^*$-members of a finite standard $A$ are standard, hence any $^*$-member of
$\bN$ which is not standard is automatically bigger than all standard
numbers.}
and by Transfer they have reciprocals in $^*\bR$, referred to as
``infinitesimal''.
\end{enumerate}

The ``philosophy'' here is to develop the theory from these properties, paying
less attention to the way such $^*V$ is constructed (or proved to exist).
The usual construction is by an ultrapower or a modification of it.
Anyway, one has to have a set $V$ to make the construction, hence only
subsets of it will have \NS\ members.
\pa

We however, choose to put the emphasis on the prescription to define the set
(or class) $\U A$ for any set (or class) $A$ -- the relation $\in$ is viewed
as just one instance of this.%
\footnote{If one wishes to be axiomatic, one may define a \NA\ in general
as a syntactic procedure to pass from a set-theoretic formula defining
$x\in A$ to a set-theoretic formula defining $x\,{}\T\in A$ (i.e.\
$x\,{}\in{}\T A$), so that the required properties are theorems.}
\pa

\section{Ultrapowers}

Let $S$ be a fixed set of indices and $\cU$ a fixed ultrafilter in $S$. For any set
$A$ we have the ultrapower $A_\cU:=A^S/\cU$, defined as the power $A^S$ modulo
the equivalence relation: equality modulo $\cU$.
\pa

We have here a covariant functor $A\mapsto A_\cU$ (which depends on $S$ and $\cU$)
from the category of sets to itself. Also, for finite $I$, $(A^I)_\cU$ is canonically
identified with $(A_\cU)^I$, and if $A'\st A$ then $(A')_\cU$ is canonically
identified with a subset of $A_\cU$, and finite unions, finite intersections
and complements carry over.
\pa

Thus one may view $A_\cU$ as the set $\U A$ of the \NS\ elements of $A$, and
the Transfer Principle is readily seen to hold w.r.t.\ first-order predicate calculus
expressions.
\pa

Note that we always have the set (or class) $\U A$ of \NS\ elements for any set (or
class) $A$. The membership relation goes as just one case. It endows any \NS\
element (=``set'') $^*x$ with the set $\al$ of its $^*$-(\NS)-members, which
by Transfer (of the extensionality of sets) determines $^*x$. The $\al$'s
obtained that way are the {\em internal} sets (of \NS\ elements). One must
emphasize that ``standard'' or ``\NS'' is here just a ``role'' that usual
mathematical (i.e.\ set-theory) entities play w.r.t.\ the construction. Sets
of \NS\ elements will naturally have their own \NS\ elements, and that can be
used to advantage.
\pa

\section{Introduction to the Proposed Construction}

In our ``philosophy'', we insist on viewing \NA\ as something totally ``part
of ordinary mathematics'', where one extends any set by adjoining all (or some)
possible limits, this in close analogy in spirit to extending a field to its
algebraic closure by adjoining all possible roots of algebraic equations.
In the latter case one would like to simply define the adjoined elements as
labeled by the irreducible polynomial they satisfy, but that is hampered by
such facts as the polynomial whose roots are the sums of the roots of two
irreducible polynomials being not irreducible.%
\footnote{Which may be viewed as the ``reason'' for the existence of several
conjugate roots: the polynomial $d(t)$ whose roots are the differences of the
roots of the irreducible $f(t)$ has a factor $t$ but also other irreducible
factors, testifying to differences $\ne0$.}
The algebraic closure thus contains conjugate elements with the same
irreducible polynomial and for any such elements there is an automorphism of
the algebraic closure exchanging them.
\pa

Similarly in our case we would like to label the adjoined ``\NS\ elements'' by
the ultrafilters on the original set, but that would not work because the
Cartesian product of two ultrafilters is in general not an ultrafilter. We too
will have a construction where there will be ``conjugate'' \NS\ elements with
the same ultrafilter and any such conjugate elements will be exchangeable by an
automorphism, provided one has enough General Continuum Hypothesis (GCH), and
then, in the same spirit as with algebraic closures, one has uniqueness (for
fixed ``basis'' $B$ -- see below) up to isomorphism, and also universality and
homogeneity.
\pa

We focus on a structure of a set $E$ being a ``\lr'' of some other, ``basis''
set $B$, a notion defined below. In fact, a \lr\ is a structure obtained in a
set $E$ whenever it serves, in any reasonable way, as a set of \NS\ elements
of $B$.
\pa

We shall construct the \NS\ elements of any set or class using
a \lr\ -- a set of \NS\ elements -- on a fixed set $B$ (i.e.\ a fixed
cardinality) and use that ``template'' to define the \NS\ elements of any set
or class $X$ (including \NS\ elements in sets of \NS\ elements!) by, in
essence, ``carrying'' them over from $B$ in some analogy with the way tangent
vectors may be defined in any smooth manifold by ``carrying'' them over from
$\bR^n$. In fact, we can acheive that goal by defining $^*X$ as a ``cylindrical''
ultrapower of $X$ w.r.t.\ to the ``cylindrical'' ultrafilter given by the \lr\
on $B$ (see below).
\pa

The cardinality of the appropriate \lr\ (set of \NS\ elements) which we shall
have on $B$ will be just $2^{2^{|B|}}$ -- the same as that of the set of
ultrafilters in $B$. For many purposes it will be convenient to choose $B$
either countable or with cardinality of the continuum, and different choices
of $B$ can be ``merged''.
\pa

We shall have the usual features: Transfer Principle, Saturation and also what
we call Confinement -- every \NS\ element belongs to some standard set of a
restricted cardinality, where the (conflicting) cardinality restrictions for
Saturation and Confinement that are obtained depend on the cardinality of $B$.
\pa

\section{\Lr s}

We say that a set $E$ has the structure of \textbf{\lr} over a set%
\footnote{In the sequel we usually assume $B$ infinite.}
(basis) $B$ if an ultrafilter $\cL$ in the cylinder Boolean algebra of $B^E$
is given. The members of the cylinder Boolean algebra are the subsets of $B^E$
that depend only on a finite number of coordinates. (In case $E$ is finite
this is just the power set of $B^E$.) What $\cL$ does is deciding, for any
(finite) family $\eta\in E^I$ indexed by a finite $I$ (pushing $\cL$, using
the map $B^E\to B^I$ induced by $\eta:I\to E$, to
an ultrafilter in $B^I$) whether relations -- subsets of $B^I$ -- hold or
not, thus transferring such $I$-relations in $B$ (= subsets of $B^I$) into 
$I$-relations in $E$, as a \NS\ setting should.
\pa

To put it otherwise: for an $n$-tuple $(e_1,\ldots,e_n)$ of elements of $E$, a
relation $R\subset B^n$\,\,$^*$-holds for $(e_1,\ldots,e_n)$ iff for $\psi\in B^E$,
$(\psi_{e_1},\ldots,\psi_{e_n})\in R$ holds modulo $\cL$.
\pa

In fact, giving $\cL$ is equivalent to giving a natural mapping between the
(contravariant) functors on sets $I$: $E^I$ and $\LP B^I\RP^*$ where the
latter denotes the set of ultrafilters in the cylinder Boolean algebra. (Note
that these functors and the natural mapping are defined for any $I$, yet they
are determined by giving them for the finite $I$.) In particular, the
``cylindrical'' ultrafilter $\cL$ in $B^E$ can be recovered as the image in
$\LP B^E\RP^*$ of the identity family in $E^E$.
\pa

But note that whenever a set $E$ serves as the set of \NS\ element of a set
$B$, with the Transfer Principle holding at least for the Propositional
Calculus operations, such a natural mapping between the above functors arises
-- a family of non-standard elements (i.e.\ an element of $E^I$ for finite $I$)
corresponds to an ultrafilter in $B^I$ consisting of the ``standard''
relations it satisfies. Thus the structure of $E$ a \lr\ over $B$ captures an
essence of the notion of a ``set of \NS\ elements of $B$''.
\pa

\begin{Rmk}
There can be sub-\lr s in two ways (or both combined):
Firstly, any subset $E'\st E$ has the structure of a \lr\ over $B$,
and secondly for any $B'\st B$ the set $\U B'$ of the $e\in E$ that satisfy
the transfer of `` belongs to $B'$ '' form a \lr\ over $B'$. Instead of
talking about sub-\lr s one may talk about embeddings. Isomorphisms and
automorphisms of \lr s are defined as expected (as those which respect the
ultrafilter $\cL$).
\end{Rmk}
\pa

\subsection {\Lr s and Ultrapowers}
\label{S:lrup}
If $E$ has the structure of \lr\ over $B$, any $e\in E$ defines an evaluation
(or coordinate) map $B^E\to B$. Thus $E$ is identified with a set of mappings
$B^E\to B$. Also, in $B^E$ we have the ultrafilter $\cL$ on the cylinder
Boolean algebra. $E$ is thus mapped into the ultrapower $B^{B^E}/\cU$, where
the ultrafilter $\cU$ is some extension of $\cL$ to all the subsets of $B^E$.
And one easily finds that the notion of a family $\eta\in E^I$ ($I$ finite)
$^*$-satisfying a relation $R\st B^I$ is the same for this ultrapower as for
the original \lr\ $E$ (whatever ultrafilter extension $\cU$ of $\cL$ to all
subsets of $B^E$ one takes).
\pa

Conversely, if $E$ is a subset of an ultrapower $B^S/\cU$ w.r.t.\ an
ultrafilter $\cU$ in $S$, then the notion of a family $\eta:I\to E$ ($I$
finite) $^*$-satisfying a relation $R\st B^I$ is defined, thus $E$ is made
into a \lr\ over $B$. To obtain the ultrafilter $\cL$ (on the cylinder
Boolean algebra of $B^E$) of this \lr, lift $E$ to a subset of $B^S$, which
makes a map $E\times S\to B$, thus a map $S\to B^E$, which pushes $\cU$ to an
ultrafilter in $B^E$. Its restriction to the cylinder Boolean algebra does
not depend on the choice of the lifting, and will be $\cL$.
\pa

Thus, since every reasonable way to endow a set $B$ with the set of its \NS\
elements $\U B$ involves a \lr, and \lr s boil down to ultrapowers, one
concludes that, in some sense, any reasonable \NA\ on a set is, in fact, defined
by a \lr\ and by an ultrapower.
\pa

\subsection{Exact \Lr s}
\label{S:exact}
For any \lr, transfer of relations will commute with Propositional Calculus
operations, but in general not necessarily with quantifiers, i.e.\ with
projections $B^{I\cup\{i\}}\to B^I$ ($i\notin I$).
\pa

This, however, will be guaranteed if we impose the requirement that
for any $\xi:I\to J$\,\,($I,J$ finite) the diagram expressing
the naturality of the above mapping between the functors is ``exact'', in the
sense that a member of $E^I$ and a member of $\LP B^J\RP^*$ which map to the
same member of $\LP B^I\RP^*$ both come from some same member of $E^J$. For
this to hold it suffices that it holds for inclusions ``adding one element''
$I\to I\cup\{i\}$ (guaranteeing that the transfer of relations commutes with
projections $B^{I\cup\{i\}}\to B^I$) and for the map $\{1,2\}\to\{1\}$
(which guarantees that the relation of equality will be transferred to the
relation of equality -- whenever this is the case for a \lr\ it is called
\textbf{separated}).
\pa

If this exactness requirement is satisfied (for all finite $I,J$), we say that
the \lr\ is \textbf{exact} or is a \textbf{\NA\ \lr}. Since for $I=\es$ both $E^\es$
and $\LP B^\es\RP^*$ are singletons with subsets naturally viewed as the
truth-values, which the natural mapping always preserves, the above commuting
of transfer of relations with logical operations implies the {\em Transfer
Principle}: any true sentence made as a first-order logical combination of
relations on $B$ will turn, by transferring each of the argument relations,
into a sentence true for $E$.
\pa

In particular, applying the exactness of the diagram for the map $\es\to\{1\}$
we have, for exact \lr s, that for every ultrafilter $\cU$ in $B$ there is an
element $e\in E$ whose ultrafilter, i.e.\ whose image in the natural mapping
$E=E^{\{1\}}\to \LP B^{\{1\}}\RP^*$, is $\cU$, that is, $e$ ``belongs'' to all
members of $\cU$ (and of course does not belong to non-members of $\cU$, which
are the complements of members).%
\footnote{This does not hold, in general, for \NA\ constructed by an ultrapower
(if the power in the ultrapower is countable, a countable set will have only
$2^{\aleph_0}$ \NS\ members but has $2^{2^{\aleph_0}}$ ultrafilters.)}
If $\cU$ is fixed (principal) -- is the family of all subsets containing
some $a\in B$ -- than $e$ is unique and is the element of $E$ to be identified
with $a$. If $\cU$ is free (non-principal), however, one proves that $e$ is
never unique and we have ``conjugate'' \NS\ elements with the same ultrafilter
which, as we shall see below, will be exchangeable by an automorphism with
favorable choice of the \lr.
\pa

But \lr s may have more: they may have also some saturation. If $m$ is a
cardinal, let us say that a \lr\ $E$ is \textbf{$m$-exact} if the above
exactness holds for the inclusion $I\to I\cup\{i\}$ if $|I|<m$. This means:
for every family $\eta:I\to E$ of elements of $E$, and any ultrafilter $\cU$
in the cylinders of $B^{I\cup\{i\}}$ which projects on $B^I$ to the
ultrafilter pulled-back from $\cL$ by $\eta$, one can find an element of
$E$ as the image of $i$ so that the pull-back of $\cL$ by the extended
$\eta$ will be $\cU$. By a transfinite process such $m$-exactness will
guarantee exactness for any inclusion $I\to M$ for $|I|<m,|M|\le m$.
In particular a \lr\ is exact iff it is separated and $\aleph_0$-exact.
It is easily seen that $m$-exactness implies $m$-saturation (see \S\ref{SCC}),
i.e.\ the property that any family of cardinality $<m$ of internal subsets of $E^n$
(i.e.\ sets of the form $(\U R)[e]$ where $R\st B^{n+1}$ and $e\in E$) has
non-empty intersection provided any finite subfamily has non-empty intersection.
Conversely, $m$-saturation implies $m_1$-exactness if $m'<m_1\Rightarrow 2^{|B|}m'<m$
(since, for infinite $B$, the cardinality of the cylinder Boolean algebra of
$B^{m'}$ is $2^{|B|}m'$). In particlar, for $m=\LP 2^{|B|}\RP^+$, $m$-exactness is
equivalent to $m$-saturation.
\pa

A separated \lr\ $E$ which is $|E|$-exact is {\em universal} in the sense that
any embedding in $E$ of a sub-\lr\ of cardinality $<|E|$ of a \lr\ $M$ (with basis
$B$) of cardinality $\le|E|$ can be extended to an embedding of all $M$, and using
a transfinite back-and-forth construction we conclude that it is {\em
homogeneous} -- any isomorphism between two sub-\lr s of cardinality $<|E|$ in
two separated \lr s $E$ with the same cardinality $|E|$ which are
$|E|$-exact can be extended to an isomorphism between the $E$'s, in
particular any two such \lr s $E$ are isomorphic.
\pa

\begin{Rmk}\label{UP}
(Combining \S\ref{S:exact} and \S\ref{S:lrup}.)
Let $E$ be an exact \lr\ over $B$.
By \S\ref{S:lrup} one constructs an ultrapower $B^{B^E}/\cL$. (By abuse of
language, we write $\cL$ for an ultrafilter on all subsets of $B^E$ extending
it, since everything will depend only on $\cL$.) $E$ embeds into the ultrapower
by mapping each $e\in E$ to the coordinate map $B^E\to B$. The fact that this
is 1-1 follows from $E$ being separated.
\pa

It might seem that the whole ultrapower must be much larger than the image of
$E$. Yet, consider the fact that $E$ is exact, and the resulting Transfer
Principle:
\pa

If $f:B^E\to B$ depends on a finite set $I\st E$ of coordinates, it is
essentially a function $f:B^I\to B$. By Transfer, one may plug $I$ into $\U f$,
and the element $e\in E$ so obtained is easily seen to be equal to $f$ modulo
$\cL$. Thus {\em any ``family'' $f:B^E\to B$ which depends only on a finite
set of coordinates is, in the ultrapower, equal modulo $\cL$ to an element of
the embedded $E$}.
\end{Rmk}

\section{Regular \Lr s; Furnishing any Set with \NNS\ Elements}

For fixed $B$ (whose cardinality $|B|$ we shall call the \textbf{base})
we can construct a separated \lr\ of cardinality $2^{2^{|B|}}$ which is
$\LP{2^{|B|}}\RP^+$-saturated, i.e.\ $\LP{2^{|B|}}\RP^+$-exact,
just by a straightforward transfinite process
indexed by the least ordinal of cardinality $\LP{2^{|B|}}\RP^+$, in each step
adjoining to $E$ all needed elements for the saturation (and extending $\cL$
correspondingly), and then taking quotient by the relation among $e,e'\in E$:
``the set $\{\xi\,|\,\xi_e=\xi_{e'}\}\st B^E$ is in $\cL$''. We call any
separated \lr\ over $B$ with this cardinality and saturation a \textbf{regular
\lr\ with basis $B$}. (Note that it was not more difficult to have any saturation
than to have ${\aleph_0}^+$-saturation.) This does not require any Generalized
Continuum Hypothesis. But if we have that, namely that
$2^{2^{|B|}}=\LP{2^{|B|}}\RP^+$ (say, if we work only in the Constructible
Universe), then a regular \lr\ $E$ will have $|E|$-exactness, hence is
universal and homogeneous as above and is unique up to isomorphism.
\pa

Any exact \lr\ $E$ over a basis $B$ induces \NS\ elements in any set or class.
They are defined as follows and will enjoy $|B|^+$-confinement where
\textbf{$m$-confinement} ($m$ a cardinality) means that any \NS\ element
$^*$-belongs to some standard set of cardinality $<m$.
(In fact, as easily seen, these will be, ``up to isomorphism'', the \NS\ elements,
that belong to standard sets of cardinality $\le|B|$, in any \NA\ for which
$E$ is the \lr\ of \NS\ elements of $B$.)
\pa

If we are to have $|B|^+$-confinement, all \NS\ elements belong to sets of
cardinality $\le|B|$, so if we know the \NS\ elements of such sets (and, by
the Transfer Principle, inclusions etc.\ among such sets will induce
inclusions of the sets of \NS\ members etc.) we will have \NS\ members of any
set or class as a direct limit.
\pa

One may think of endowing sets $X$ of cardinality $\le|B|$ with
\NS\ elements by ``carrying'' them over from $B$ in some analogy with the way
tangent vectors may be defined in any smooth manifold by ``carrying'' them
over from $\bR^n$ -- just define a \NS\ element of $X$ as a mapping which
corresponds to each bijection $j:X\simeq B'\st B$ an element $e_j$ of $E$
which $^*$belongs to $B'$, so that for different bijections $j$ the $e_j$ are
connected by the transformation that connects the $j$'s. The Transfer
Principle holding in the \lr\ $E$ will ensure that everything fits.
\pa

But the same effect can be achieved more neatly (see Remark \ref{UP}) by
defining the \NS\ elements of any set $X$ as the members of a ``cylindrical''
ultrapower $X^{B^E}/\cL$ where one takes the set of all maps $B^E\to X$
{\em that depend only on a finite number of coordinates} and factors by
equivalence modulo $\cL$.
\pa

In this way, for every set (or class) $A$ we will have the ``mirror'' set
(or class) $\U A$ of the non standard elements of $A$ (all that in our usual
set theory). We will have the Transfer Principle (i.e.\ commuting of
first-order logical operations with the transfers $A\to{\U A}$) in general
($|B|^+$-confinement easily reduced everything to working within a set of
cardinality $|B|$). In a similar way, we will have $m$-saturation in general
if we have it in $B$ {\em provided $m\le|B|^+$}. When $E$ is a regular \lr\
with basis $B$ and one is working inside a standard set of cardinality $\le|B|$
(like when one works inside $B$) one has saturation $\LP{2^{|B|}}\RP^+$.
\pa

A \NS\ element has its degree of confinement -- how small the cardinality is
of a standard set to which it belongs. Note that anyway, for a cardinal $m$,
when one works with standard sets of cardinality $\ge m^+$ one cannot have
confinement $\le m^+$ with saturation $\ge \LP m^+\RP^+$ -- with
$\LP m^+\RP^+$-saturation there will be a \NS\ $^*x$ which does not belong
to any standard set of cardinality $\le m$. Indeed, well-order a standard set of
cardinality $m^+$ by the set of ordinals with cardinality $\le m$ and choose
(by saturation) a $^*x$ greater then all standard elements. Thus confinement
and saturation are somehow conflicting requirements, and for various needs one
may take different bases (see below).
\pa

We shall fix a regular \lr\ over a set $B$ (referred to as the
{\EM template}) with cardinality $|B|=b$ (called the {\EM base}) to endow any
set (or class) with \NS\ elements. Base $b$ will give confinement
$b^+$, saturation $b^+$ in general but saturation
$\LP2^b\RP^+$ if one works inside a standard set of cardinality $\le b$. The
number of \NS\ elements of any infinite standard set of cardinality
$\le2^{2^b}$ will be $2^{2^b}$. For infinite $b'<b$, there will be among them
$2^{2^b}$ \NS\ elements with ``better'' confinement $b'$ (see the cardinality
computations in the next \S).
\pa

For two regular \lr s with different bases $b_1<b$, the base $b_1$ \lr\ will have
cardinality $2^{2^{b_1}}$ while the base $b$ \lr\ will have saturation, equivalently
exactness, $\LP2^b\RP^+$. Suppose, e.g.\ that $b\ge2^{b_1}$ (say $b_1=\aleph_0$
and $b=$ continuum) then by saturation we can embed the template $b_1$-\lr\
into the template $b$-\lr\ and if we have the case $\LP2^b\RP^+=2^{2^b}$ of
GCH then any two such embeddings are exchanged by an automorphism of the
$b$-\lr. Thus in this case the $b_1$-\NS\ elements can be viewed as identified
with part of the $b$-\NS\ elements (indeed part of those with confinement
$\le {b_1}^+$). This identification is not unique but with GCH as above all
identifications are ``equivalent''. In this way working with $b_1$-\NS\
elements (in any set or class) merges smoothly with working with $b$-\NS\
elements.

\section{More about Saturation, Some Calculations of Cardinalities}
\label{SCC}

As mentioned above, in \NA, an {\EM internal set} of \NS\ elements is
one obtained as the set $\al$ of \NS\ $^*$-members of some \NS\
($^*$-set) $^*x$. $\al$ determines $^*x$ (by transfer of extensionality),
and may serve as a substitute when talking about $^*$-sets: in transferring a
statement about sets translate: set (of standard elements) $\to$ internal set
(of \NS\ elements). One easily shows that equivalently an internal set
is one obtained as a section (parallel to one of the axes) of $\U R$, for a
standard set $R$ of ordered pairs, through a \NS\ element.
In general a set of \NS\ elements is not internal: indeed, since,
by Transfer, every bounded non-empty internal subset of $^*\bN$ has a maximal
element, the set of $^*$-natural-numbers that are identified with the standard
ones is not internal. A non-internal subset of a set $\T A$ (for some standard $A$)
is called {\EM external}.
\pa

The assumption of $m$-exactness says that if $X$ is a standard set,
$\U X$ the \lr\ over $X$ of all \NS\ members of $X$,\,$|I|<m$,
$i_0\notin I$, $\eta:I\to{}\U X$ and $\cU\in\LP X^{I\cup\{i_0\}}\RP^*$
such that the projection of $\cU$ on $X^I$ coincides with the pullback of
$\cL$ by $\eta$, then $\exists$ a $\U x\in{}\U X$ so that if we extend $\eta$
by mapping $i_0\to{}\U x$ then the pullback will be $\cU$.
\pa

Now, the members of $\cU$ depend each on a finite number of coordinates,
and as such are standard relations. When we substitute the
$\eta(i)$'s we get internal properties (=sets) of the sought-for $\eta(i_0)$.
These properties are not empty, since the $(\eta(i))_i$ $^*$-belong to the
projections of the members of $\cU$ (here we use Transfer). Therefore, since
$\cU$ is an ultrafilter, any finite collection of these properties intersect.
The required existence of a suitable $\U x=\eta(i_0)$ means that all intersect.
Elaborating this, one connect $m$-exactness with $m$-saturation, i.e.\ the property
that:
\pa
\begin{tabular}{ll}
&For any family of less than $m$ internal sets (= internal properties),\\
&if every finite subfamily has non-empty intersection,\\
&then the whole family has non-empty intersection.
\end{tabular}
\pa

As an example for the use of saturation, let us prove
\begin{Prop}
Let us work with a \NA\ as above with base $b$. Then
\begin{itemize}
\item[(i)]
Let $B$ be a set of cardinality $b$. Then for every internal subset
$\be\st\T B$ which contains all standard elements of $B$, and for every
ultrafilter $\cU$ on $B$, $\exists\;\U x\in\be$ whose ultrafilter is $\cU$.
Consequently, $|\be|=2^{2^b}$ (recall that $|\U{}B|=2^{2^b}$).
\item[(ii)]
If $\al$ is internal, then either $\al$ is finite or $\exists$ an internal
one-one
mapping from an internal $\be$ as in (i) to $\al$, hence $|\al|\ge2^{2^b}$.
\item[(iii)]
If $\al$ is $^*$-finite (i.e.\ a $^*$-member of a standard family of finite
sets) then $|\al|$ is finite or $2^{2^b}$.
\end{itemize}
\end{Prop}
\begin{Prf}
(i) We use $\LP2^b\RP^+$-saturation (which holds when working inside $B$):
for any finite family of members of $\cU$ $\exists\;\U x\in\be$ belonging
to all of them (indeed a standard one), and by saturation we are done.

(ii) Here we use $b^+$-saturation: if $\al$ is not finite, then for every
finite set $F$ of standard members of $B$\,\,\,$\exists$ an internal set
$\be_F\st{}\U B$ containing all members of $F$ and an internal one-one mapping
$\be_F\to\al$ (take $\be_F=F$), hence by saturation $\exists$ a $\be$ suitable
for all standard members of $B$.

(iii) Suppose $\al$ not finite. By (ii) $|\al|\ge2^{2^b}$. On the other hand,
by $b^+$-confinement $\al$ $^*$belongs to a standard set $S$ of cardinality
$\le b$ whose members are finite sets. Then $|\cup S|\le b$, thus
$|{}^*(\cup S)|\le2^{2^b}$, while $\al\st{}^*(\cup S)$.
\end{Prf}

\end{document}